\documentclass[12pt,draftcls,onecolumn]{IEEEtran}

\IEEEoverridecommandlockouts

\usepackage{latexsym}
\usepackage{graphicx}
\usepackage{epsfig}
\usepackage{subfigure}
\usepackage{array}
\usepackage{amsmath}
\usepackage{amssymb}
\usepackage{color, soul}
\usepackage{amsthm}
\usepackage{enumerate}
\usepackage{algpseudocode}
\usepackage{algorithm}

\allowdisplaybreaks

\theoremstyle{plain} 

\theoremstyle{definition}

\newcommand{\bI}{{\bf I}}
\newcommand{\bU}{{\bf U}}
\newcommand{\bV}{{\bf V}}

\newcommand{\bA}{{\bf A}}
\newcommand{\bB}{{\bf B}}

\newcommand{\bp} {\begin{proof}}
\newcommand{\ep} {\end{proof}}

\newcommand{{\Rb}} {\right)}

\newcommand{{\Rf}} {\right\}}

\def\bal#1\eal{\begin{align}#1\end{align}}

\begin{document}

\title{On Singular Value Inequalities for the Sum of Two Matrices}

\author{Sergey Loyka
%
\thanks{S. Loyka is with the School of Electrical Engineering and Computer Science, University of Ottawa, Ontario,
Canada, K1N 6N5, e-mail: sergey.loyka@ieee.org.}
}
%
\maketitle
%

\vspace*{-1.5\baselineskip}

\begin{abstract}
A counter-example to lower bounds for the singular values of the sum of two matrices in [1] and [2] is given. Correct forms of the bounds are pointed out.
\end{abstract}



%
\vspace*{1\baselineskip}

The following lower bound for singular values of the sum of two $n \times n$ matrices $\bA$ and $\bB$ is reported in \cite{Marshall'11} (see G.3.b on p. 334):
\bal
\sum_{i=1}^{k} \sigma_i(\bA+\bB) \ge \sum_{i=1}^{k} \sigma_i(\bA) - \sum_{i=1}^{k} \sigma_{n-i+1}(\bB),\ 1 \le k \le n,
\eal
where singular values $\sigma_i$ are in decreasing order, $\sigma_1 \ge \sigma_2 \ge ...$. This lower bound is also claimed in \cite{Zhang'11} (see Problem 4 in Sec. 10.5, p. 361). For $k=1$, this reduces to
\bal
\label{eq.sigma_1}
\sigma_1(\bA+\bB) \ge \sigma_1(\bA) - \sigma_{n}(\bB)
\eal
To see that this does not hold in general, consider the following example: $\bA=diag\{1,0\}, \bB=diag\{-1,0\}$, which results in
\bal
\label{eq.sigma_1}
0= \sigma_1(\bA+\bB) \ge \sigma_1(\bA) - \sigma_{n}(\bB) = 1
\eal
This example also disproves the lower bound given in Theorem 8.13 of \cite{Zhang'11}, namely
\bal
\sigma_i(\bA+\bB) \ge \sigma_i(\bA) + \sigma_{n}(\bB)
\eal
The correct form of this bound is
\bal
\sigma_i(\bA+\bB) \ge \sigma_i(\bA) - \sigma_{1}(\bB)
\eal
which follows from e.g. Theorem 3.3.16(c) in \cite{Horn'91}.

An inspection of the "proof" in \cite{Marshall'11} reveals the following incorrect step
\bal
\min \mathrm{Re}\ \mathrm{tr} \bU\bB\bV^+ = -\sum_{i=1}^k \sigma_{n-i+1}(\bB)
\eal
where $\min$ is over all semi-unitary matrices $\bU, \bV$ such that $\bU\bU^+=\bV\bV^+=\bI_k$ and $\bI_k$ is $k\times k$ identity matrix.
The correct form of this step is as follows:
\bal
\label{eq.step.correct}
\min \mathrm{Re}\ \mathrm{tr} \bU\bB\bV^+ = - \max \mathrm{Re}\ \mathrm{tr} \bU(-\bB)\bV^+ = -\sum_{i=1}^k \sigma_{i}(-\bB) = -\sum_{i=1}^k \sigma_{i}(\bB)
\eal
where 2nd equality is due to e.g. Theorem 3.4.1 in \cite{Horn'91}, and the last equality is due to $\sigma_{i}(-\bB) = \sigma_{i}(\bB)$, which results in the following lower bound:
\bal
\label{eq.bound.correct}
\sum_{i=1}^{k} \sigma_i(\bA+\bB) \ge \sum_{i=1}^{k} \sigma_i(\bA) - \sum_{i=1}^{k} \sigma_{i}(\bB)
\eal

This lower bound is not the tightest one in general. A tighter lower bound can be established using e.g. Theorem 3.4.5. in \cite{Horn'91}:
\bal
\label{eq.T10.24}
\sum_{i=1}^{k} \sigma_i(\bA-\bB) \ge \sum_{i=1}^{k} s_{[i]}
\eal
where $s_i = |\sigma_{i}(\bA) - \sigma_{i}(\bB)|$ and $s_{[i]}$ are in decreasing order (note that while $\sigma_{i}$ are in decreasing order, $s_i$ do not need to be so), from which it follows that
\bal
\label{eq.better} \notag
\sum_{i=1}^{k} \sigma_i(\bA+\bB) &\ge \sum_{i=1}^{k} s_{[i]}\\ \notag
&= \max_{1\le i_1 <...<i_k\le n} \sum_{l=1}^{k} |\sigma_{i_l}(\bA) - \sigma_{i_l}(\bB)|\\ \notag
&\ge \max_{1\le i_1 <...<i_k\le n} \sum_{l=1}^{k} (\sigma_{i_l}(\bA) - \sigma_{i_l}(\bB))\\
&\ge \sum_{i=1}^{k} (\sigma_{i}(\bA) - \sigma_{i}(\bB))
\eal
where 1st inequality follows from \eqref{eq.T10.24} via the substitution $\bB \rightarrow -\bB$.

It should be noted that \eqref{eq.step.correct}-\eqref{eq.better} hold for rectangular matrices as well.


\end{document}